\documentclass[12pt]{amsart}

\usepackage{amssymb}

\newcommand{\La}{\langle}
\newcommand{\Ra}{\rangle}

\newcommand{\ran}{\operatorname{Ran}}

\newcommand{\wt}{\widetilde}

\newcommand{\fdot}{\,\cdot\,}

\newtheorem{thm}{Theorem}[section]
\newtheorem{lm}[thm]{Lemma}

\theoremstyle{remark}

\newtheorem{rem}[thm]{Remark}

\newtheorem*{rem*}{Remark}

\theoremstyle{definition}
\newtheorem{dfn}[thm]{Definition}
\newtheorem*{dfn*}{Definition}

\newcommand{\ci}[1]{_{ {}_{\scriptstyle #1}}}

\newcommand{\ti}[1]{_{\scriptstyle \text{\rm #1}}}

\newcommand{\cD}{\mathcal{D}}
\newcommand{\cS}{\mathcal{S}}

\newcommand{\R}{\mathbb{R}}
\newcommand{\Z}{\mathbb{Z}}
\numberwithin{equation}{section}

\makeatletter
\def\multilimits@{\bgroup
   \Let@
   \restore@math@cr
   \default@tag
  \baselineskip\fontdimen10 \scriptfont\tw@
  \advance\baselineskip\fontdimen12 \scriptfont\tw@
  \lineskip\thr@@\fontdimen8 \scriptfont\thr@@
  \lineskiplimit\lineskip
  \vbox\bgroup\ialign\bgroup\hfil$\m@th\scriptstyle{##}$\hfil\crcr}
\def\Sb{_\multilimits@}
\def\Sp{^\multilimits@}
\def\endSb{\crcr\egroup\egroup\egroup}

\makeatother

\newcommand{\cz}{Calder\'{o}n--Zygmund\ }

\normalbaselines

\begin{document}

\title[Inequalities for Haar multipliers]{Two weight inequalities for 
individual Haar multipliers and other
well localized operators }
\author{F.~Nazarov, S.~Treil and A.~Volberg}
\thanks{This paper is based upon work supported by the National Science Foundation under Grant  DMS-0501065.
}
\date{}

\begin{abstract}
In this paper we are proving that Sawyer type condition  for boundedness 
work for the two weight estimates of individual
Haar multipliers, as well as for the Haar shift and other ``well 
localized'' operators.
\end{abstract}

\maketitle

\setcounter{section}{-1}
\section{Introduction}
\label{intro}

The main question of this paper concerns two weight estimates for singular 
integral  operators,
i.e.~the questions when an integral operator $T$ is bounded operator 
from a weighted
space $L^2(w)$ to $L^2(v)$.

One of the most interesting cases is the case when $T$ is the Hilbert
transform, $Tf(s) = \frac1\pi \int\frac{f(t)}{s-t} dt$, although the case of more general \cz operators seems to be of great interest as well. For such operators  no ``real variable'' necessary and sufficient condition is known.%
\footnote{For the Hilbert transform a necessary and sufficient condition of Helson--Szeg\"{o} type was obtained by M.~Cotlar and C.~Sadosky \cite{CS1}. Their condition was stated in the language of complex analysis: the Hilbert transform is a bounded operator from $L^2(w)\to L^2(v)$ if and only if one can find a function $h$ in the analytic Hardy class $H^1$ such that for some $C>0$ the matrix $\left(\begin{array}{cc}
Cw-v & Cw+v-h \\ 
Cw+v-\overline h & Cw-v \\ 
\end{array}\right)$ 
is positive semidefinite a.e. No generalization of this condition to other \cz operators is known. 
}

In this paper we deal with dyadic analogues of the Hilbert transform, the
so-called {\em Haar multipliers} and their generalizations. It turns out
that for such operators it is possible to find  necessary and sufficient
condition for two weight estimates.

Let us introduce one of the main examples. The {\em standard dyadic grid}
$\cD=\cD_0$ is the collection of all {\em dyadic intervals} $[2^k\cdot j,
2^k\cdot (j+1))$, $k, j\in \Z$. A general dyadic grid is an object
obtained from $\cD_0$ by a dilation and a shift.

For an interval $I\subset \R$ we define the ($L^2$-normalized) Haar
function  $h\ci I:= |I|^{-1/2}(\chi\ci{I_+}-\chi\ci{I_-})$;
here $|I|$ stands for the length of the interval $I$; $I_+$ and $I_-$
denote its right and left halves respectively, and $\chi\ci E$ denotes the
characteristic function (indicator) of the set $E$.

Given a sequence $\mathbf\alpha = \{\alpha\ci I\}\ci{I\in \cD}$ define
the Haar multiplier (a.k.a.~the martingale transform) $T^{\mathbf\alpha}$ by
$$
T^\alpha f := \sum_{I\in \cD} \alpha\ci I (f, h\ci I) h\ci I, \qquad 
f\in L^2(\R).
$$
We are interested  under what conditions on the weights $v$, $w$ ($v, 
w\ge0$, $v,
w\in L^1_{\text{loc}}$) the operator $T=T^\alpha$ extends to a bounded operator from 
$L^2(w)$ to
$L^2(v)$, i.e., when the following {\em two weight estimate}
$$
\int_\R |T^\alpha f|^2 v dx \le C\int_\R |f|^2 w dx
$$
holds.

If one assumes that $w^{-1}\in L^1_{\text{loc}}$ (this assumption is 
satisfied in
practically all interesting cases), then the above estimate means that 
the operator
$M_{v^{1/2}} T^\alpha M_{w^{-1/2}}$, where $M_\phi$ stands for the 
operator of
multiplication by $\phi$, is bounded in the non-weighted $L^2$ space.

In \cite{NTV-2weight} this question was studied for the family of Haar 
multipliers
$T^{\sigma\alpha}$, where $\alpha= \{\alpha\ci I\}\ci{I\in \cD}$ is 
fixed and
$\sigma= \{\sigma\ci I\}\ci{I\in \cD}$ is an arbitrary sequence of signs,
$\sigma\ci I =\pm1$. It was shown that the operators
$M_{v^{1/2}}T^{\sigma\alpha}M_{w^{-1/2}}$ are uniformly (over all 
possible choices
of signs
$\sigma$) bounded if and only if the operators  are uniformly bounded on 
the test
functions  $w^{-1/2}\chi\ci I$, $I\in \cD$, and their adjoints
$M_{w^{-1/2}}T^{\sigma\alpha}M_{v^{1/2}}$ are uniformly bounded on the
test
functions  $v^{1/2}\chi\ci I$, $I\in \cD$.

Conditions where the boundedness of an operator follows from its 
boundedness on
special test functions as above are called ``Sawyer type conditions'',
after
E.~Sawyer, who proved in \cite{Saw1} that such condition is necessary and
sufficient for two weight estimates of the maximal function (only one 
condition is
needed there). After that there were quite a few results that Sawyer type
conditions are sufficient for the boundedness of  (clearly they are 
necessary) of
different classes of integral operators with {\em positive kernels}.

The above mentioned result about Haar multipliers was a first (and up 
until now a
unique)  result giving necessary and sufficient condition for the boundedness of operators with alternating kernels. The fact 
that helped us a
lot   was that we were dealing with a {\em family} of operators, and
that allowed us to reduce the problem to estimates of an operator with
non-negative kernel.

At that moment it seemed impossible to get a sharp result about two weight
estimates for an {\em individual} Haar multiplier, at least our method 
did not
allow us to do that.

So, the main result of the paper looks quite surprising: the Sawyer type 
condition
works for two weight estimates of an individual Haar multiplier (as well 
as for a
wider class of so-called  ``band operators'')!

Before giving a formal definition, let us present one more important 
example of an ``band operator'',  the so called \emph{Haar 
Shift} $\cS$
$$
\cS f := \sum_{I\in \cD} (f, h\ci I) [h\ci{I_+}-h\ci{I_-}] .
$$
This operator is interesting, in particular, because if we average it 
over all dyadic grids (translated and rescaled) we get the Hilbert 
transform (up to a multiplicative constant), see \cite{Pet}.

\section{Two weight estimates of band operators: formal definitions and main results}

\subsection{Main definitions}
Let $\cD$ be the standard dyadic lattice in $\R^N$.
Let us recall the definition: for
each $k\in\Z$ we consider the cube $[0, 2^k)^N$ and all its translations
by elements of $\R^N$ with coordinates of form $j\cdot 2^k$,
$j\in\Z$, then  take union over all $k$.
Each dyadic cube has $2^N$ ``sons'', and each cube has a unique
``parent''.

In this paper by a cube we will always mean a dyadic cube.

For a cube $Q$ let $\ell(Q)$ denote its size, i.e. the length of its 
side. And for a cube $Q$ let $Q^{(k)}$ be the $k$th grandparent of the 
cube $Q$, i.e.~the cube $R$ of size $2^k\ell(Q)$ containing $Q$. 

A (non-weighted) Haar function on a cube $Q$ is a function $h\ci Q$ supported on $Q$, constant on the ``children'' of $Q$ and orthogonal to constants, 
$\int_Q h\ci Q dx =0$. The set of all Haar function on a cube $Q\subset \R^N$ is a vector space of dimension $2^N-1$. 

With a dyadic lattice $\cD$ in $\R^N$ one can naturally associate a (non-oriented) $2^N$-ary tree, where each dyadic cube is connected to its $2^N$ ``sons''. By the ``tree distance'' (or ``graph distance'') $d_\text{tree}(Q, R)$ between cubes $Q, R\in\cD$ we understand the distance on the graph, where we assign length $1$ to each wedge. 
  
And now the formal definition:
\begin{dfn}
A \emph{band operator}  on $\R^N$  is a bounded operator in $L^2(\R^d)$  
whose matrix in the Haar basis has the ``band structure'', meaning there exists $r\in \mathbb{Z}_+$ such that $(Th\ci 
Q ,h\ci R)=0$ for all Haar functions $h_Q$, $h_R$ such that $d_\text{tree}(Q,R)>r$. 
\end{dfn}

The Haar multipliers (martingale transforms) discussed in the introductions are band operators with $r=0$, and the Haar shift is a band operator with $r=1$. 

\subsection{Two weigh estimates}
As we already mentioned above, the operator $T$ can be extended to a bounded operator acting from $L^2(w)$ to $L^2(v)$ if and only if the operator $M_v^{1/2} T M_w^{-1/2}$ can be extended to a bounded operator in the (non-weighted) $L^2$. Denoting $u:=w^{-1}$ we can write the problem in more symmetric form, namely, when the operator $M_v^{1/2} T M_u^{1/2}$ can be extended to a bounded operator in the non-weighted $L^2$? Here we assume that $u$ and $v$ are $L^1_\text{loc}$ weights. This formulation is (at least formally) a bit more general then the formulation with two weight estimates, because here we do not assume that $u^{-1}\in L^1_\text{loc}$ (and we did assume above that $w^{-1} \in L^1_\text{loc}$).

Now we are ready to state one of the main results of the paper:

\begin{thm}
\label{t0.1}
Let $T$ be a band operator, and let $u, v\ge 0$, $u, v\in L^1_{\text{loc}}$. Then the operator $M_v^{1/2} T M_u^{1/2}$ extends to a bounded operator in $L^2$ if and only if for all dyadic cubes $Q$ the following two inequalities hold:
\begin{align}
\label{0.1}
\int_{\R^N} | T(\chi\ci Q u) |^2 v dx  &\le C\int_{Q} u dx ,\\
\label{0.2}
\int_{\R^N} | T^*(\chi\ci Q v) |^2 u dx  &\le C\int_{Q} v dx.
\end{align}
Moreover, the norm of the operator $M_v^{1/2} T M_u^{1/2}$ can be estimated by a constant depending only on the dimension $N$, the number $r$ in the definition of the band matrix, and the constants $C$ in \eqref{0.1}, \eqref{0.2}
\end{thm}

The above theorem says that to verify the boundedness of the operator $M_v^{1/2} T M_u^{1/2}$ it is sufficient only to check its boundedness  on the test functions $u^{1/}\chi\ci Q$ and the boundedness of its formal adjoint $M_u^{1/2} T^* M_v^{1/2}$ on the test functions $v^{1/2} \chi\ci Q$ for all dyadic cubes $Q$.

\begin{rem}
\label{r0.2}
A careful reader can ask a question here on how one can interpret the expressions $T(\chi\ci Q u)$ and $T^*(\chi\ci Q v)$. The problem is that the operator $T$ (and so its adjoint $T^*$) are  defined only on $L^2$, so if one only assumes that $u, v\in L^1_\text{loc}$, the above expressions formally are not defined. 

However, it is easy to give the meaning to the conditions \eqref{0.1}, \eqref{0.2}. First of all, $T(\chi\ci Q u)$ and $T^*(\chi\ci Q v)$ are well defined if $u, v\in L^2_\text{loc}$. Second, it is trivial, that for any weight $v'\le v$ we have $\|M_{v'}^{1/2} T M_u^{1/}\| \le \|M_{v}^{1/2} T M_u^{1/}\|$, and by taking adjoint one can conclude that for any weight $u'\le u$ we have $\|M_{v}^{1/2} T M_{u'}^{1/}\| \le \|M_{v}^{1/2} T M_u^{1/}\|$. 

Therefore, for any sequence of weights $u_n$, $u_n\in L^2(Q)$, $u_n \nearrow u$, the expression $Tu_n\chi\ci Q$ is well defined.  If, in addition, one assumes that the operator $M_v^{1/2} T M_u^{1/2}$ is bounded, then $\| v^{1/2} T u_n \chi\ci Q \|_{L^2}\le C<\infty$, and moreover $v^{1/2} T u_n \chi\ci Q$ converges to some functions in $L^2$. It is easy to see that this limit function does not depend on the choice of the sequence $u_n$, so we call this function $v^{1/2} T\chi\ci Q u$. The right side of \eqref{0.1} then can be interpreted as the $L^2$ norm of this function. The dual condition \eqref{0.2} can be interpreted similarly. 

\end{rem}

\begin{rem}
For the sufficiency part one can think even of a simpler condition to interpret
 \eqref{0.1}, \eqref{0.2}. For example, one can pick an increasing sequence of weights $u_n\in L^2_\text{loc} \nearrow u$ (for example $u_n(x) = \min\{u(x), n\}$) and interpret the condition \eqref{0.1} as the uniform estimate (independent of $n$ and $Q$)
$$
\int_{\R^N} | T(\chi\ci Q u_n) |^2 v dx  \le C\int_{Q} u dx
$$
The dual condition \eqref{0.2} is interpreted similarly by picking an increasing sequence $v_n\nearrow v$. 

Indeed, these conditions would imply the conditions \eqref{0.1}, \eqref{0.2} with $u$ and $v$ replaced by $u_n$ and $v_m$ respectively. That, by Theorem \eqref{t0.1} implies the uniform estimate $\| M_{v_n}^{1/2}T M_{u_m}^{1/2}\|\le C<\infty$ which implies the estimate $\| M_{v}^{1/2}T M_{u}^{1/2}\|\le C$.
\end{rem}

The conditions of Theorem \ref{t0.1} can be relaxed a bit by integrating only over the cubes $Q$:

\begin{thm}
\label{t0.3}
Let $T$ be a band operator  and let $u, v\ge 0$, $u, v\in L^1_\text{loc}$. Then the operator $M_v^{1/2} T M_u^{1/2}$ extends to a bounded operator in $L^2$ if and only if for all dyadic cubes $Q$ the following two conditions  hold:
\begin{enumerate}
\item For all dyadic cubes $Q$, $R$ satisfying $2^{-r}\ell(Q)\le \ell(R) \le 2^r \ell(Q)$
$$
|\langle  T \chi\ci Q , \chi\ci R  \rangle | \le C \left(\int_Q u \right)^{1/2} \left(\int_Q v \right)^{1/2} ;
$$
\item For all dyadic cubes $Q$
$$
\int_{\R^N} | T(\chi\ci Q u) |^2 v dx  \le C\int_{Q} u dx ,\quad
\int_{\R^N} | T^*(\chi\ci Q v) |^2 u dx  \le C\int_{Q} v dx.
$$
\end{enumerate}
Moreover, the norm of the operator $M_v^{1/2} T M_u^{1/2}$ can be estimated by a constant depending only on the dimension $N$, the number $r$ in the definition of the band matrix, and the constants $C$ in the above conditions (1) and (2). 
\end{thm}

In this paper $\langle \fdot, \fdot\rangle$ stands for the standard inner product in $L^2(\R^N)$, $\La f, g\Ra = \int f\overline g dx$. Note that 
$$
\La T \chi\ci Q u , \chi\ci R v \Ra = \La M_v^{1/2} T M_u^{1/2} \chi\ci Q u^{1/2} , \chi\ci R v^{1/2} \Ra, 
$$
so, as we discussed above in Remark \ref{r0.2}, this expression is well defined. 

\begin{rem}
Note, that in fact one does not even need the operator $T$ to be bounded in $L^2$. A bit more elaborate reasoning than in Remark \ref{r0.2} would allow us to interpret the conditions of Theorem \ref{t0.1} in the case when we only require  the bilinear form $\La Tf, g\Ra$ of the operator $T$ be defined on bounded compactly supported functions $f$ and $g$. We leave the details here as an exercise for the reader. 
\end{rem}

\section{Well-localized operators for general measures}

To prove theorems \ref{t0.1} and \ref{t0.3}  we prove even a bit more general results about the so-called \emph{well-localized} operators $T=T_\mu:L^2(\mu)\to L^2(\nu)$ for general measures $\mu$, $\nu$ that may have singular parts.

\subsection{Heuristics and formal definition}The idea behind the notion of well-localized operators is rather simple. If one has an integral operator $T$, 
$$
Tf(x)=\int K(x, y) f(y) dm(y), 
$$
where the integration is with respect to the Lebesgue measure $m$ in $\R^N$, one can consider (at least formally) the operator $T_\mu$ where the integration with respect to the measure $\mu$, 
$$
T_\mu f(x)=\int K(x, y) f(y) d\mu(y).
$$
Such reduction is quite common in weighted estimates. For example, two weight estimates for integral operators $T$ can be reduced to estimates of operators $T_\mu:L^2(\mu)\to L^2(\nu)$ with appropriately chosen measures $\mu$ and $\nu$. One of the advantages of such representation is that the integration in the operator and in the computation of the norm is with respect to the same measure $\mu$. 

Since a rank one operator $\La\fdot,f\Ra g$ is an integral operator with kernel $g(x) \overline{f(y)}$, one can formally represent any operator $T$ is $L^2$ as an integral operator by considering its matrix in an orthonormal basis (for example in the Haar basis) and writing it as a sum of rank one operators. And then one can replace the integration with respect to the Lebesgue measure by the integration with respect to the measure $\mu$. However this is only a formal reasoning, because, first of all, the resulting ``kernel'' of the ``integral'' operator does not need to be a function. And even if the kernel $K$ is a function, it is completely not clear how to interpret the operator $T_\mu$. 

So, instead of  writing a band operator as an integral operator (which is not always possible, for example the identity $I$ is a band operator) and then trying to interpret the operator $T_\mu$, we pick an \emph{axiomatic approach}. 

Looking at the above  formal reasoning, one can figure out  what structure the matrix of the operator $T_\mu$ should have with respect to the weighted Haar bases in $L^2(\mu) $ and $L^2(\nu)$, if one start with a band operator $T$. It is not difficult to see, that while the matrix of $T_\mu$ in the weighted Haar bases is not generally a band matrix, it has some special properties, some traces of the band structure of $T$. 

So, we took this special structure of the matrix in the weighted Haar bases as the definition of the \emph{well-localized} operators. Later in this section we will show how for a band operator $T$ one can rigorously reduce the estimates of  the operator $M_v^{1/2} T M_u^{1/2}$ to the estimates of of the appropriate  well-localized operator $T_\mu$

And now main definitions. 

First, let us define the weighted Haar system. 
For each cube $Q$ a Haar function with respect to a measure $\mu$ 
($\mu$-Haar function) on $Q$ is a function $h^\mu_Q$ supported on $Q$, constant 
on all $2^N$ ``children'' of $Q$, and such that $\int_Q h_Q^\mu \,d\mu =0$. 
Note that for a given $Q$ the set $H_Q^\mu$ of all $\mu$-Haar functions on $Q$ is a 
subspace of dimension at most $2^N-1$ (can be less since degenerate 
situations are possible). The set $H_Q$ of non-weighted Haar functions on $Q$ has exactly dimension $2^N-1$.

Let $\mu$ and $\nu$ are (regular Borel) measures in $\R^N$.
Let us be given an operator $T=T_\mu$ acting from $L^2(\mu)$ to 
$L^2(\nu)$. By given we mean that we know its bilinear form $\La T \chi\ci Q, \chi\ci R\Ra_\nu$ on characteristic functions of cubes; here $\La\,\cdot\,, \,\cdot\,\Ra_\nu$ is the inner product in $L^2(\nu)$, 
$$
\La f, g\Ra_\nu =\int f\overline g\,d\nu . 
$$
Note that the above bilinear forms also define a formal adjoint $T^*=T^*_\nu$ of $T=T_\mu$.

\begin{dfn}\label{df1.1}
Let  $T=T_\mu$ be he operator defined above,  acting (formally) from 
$L^2(\mu)$ to $L^2(\nu)$. We say that $T_\mu$ is  \emph{lower triangularly localized} if 
there exists a constant $r>0$ such that for 
all cubes $R$ and $Q$, $\ell(R)\le \ell(Q)$,  and for all $\nu$-Haar 
functions $h^\mu_R$ on $R$
$$
\La T_\mu \chi\ci Q , h_R^\nu\Ra_\nu =0
$$
if $R\not\subset Q^{(r)}$ or if $\ell(R) \le 2^{-r} \ell(Q)$ and 
$R\not\subset Q$. Here, recall, $Q^{(r)}$ is the ``grandparent'' of 
order $r$ of the cube $Q$.

And we say that the operator $T_\mu$ is \emph{well localized} if both 
$T_\mu$ and its formal adjoint $T^*_\nu$ are lower triangularly localized.
\end{dfn}

Note, that the Haar multipliers and Haar shift, discussed in the 
Introduction are well localized in the sense of the above definition, see Section \ref{s1.1} below. 
Note, that while the matrix of a Haar multiplier $T^\alpha$ in the 
non-weighted case has only one diagonal, the matrix of the weighted 
version may have infinitely many.

\subsection{From weighted estimates for band operators to estimates of well localized operators}
\label{s1.1}

First we want to show that two weight estimates for the band operators can be reduced from  the estimates for the well-localized operators. 

Let $T$ be a band operator in the non-weighted $L^2(\R^N)$ and let $u, v\ge 0$, $u, v\in L^1_\text{loc}$ be two weights. Suppose we want to know whether the operator $M_v^{1/2}T M_u^{1/2}$ is bounded.

Let us denote $d\nu = vdm$, $d\mu = u dm$, where $m$ is the Lebesgue 
measure in $\R^N$. Note that $M_{u^{1/2}}: L^2(\mu)\to L^2(m)$, $M_{v^{1/2}}: L^2(\nu)\to L^2(m)$ are isometries. So the boundedness of the 
operator $M_v^{1/2}TM_u^{1/2}$ in $L^2(m)$ is equivalent to 
the boundedness of the operator
$$
T_\mu = M_{v^{-1/2}} (M_v^{1/2} T M_u^{1/2} ) M_{u^{1/2}} = T M_{u}
$$
acting from $L^2(\mu)$ to $L^2(\nu)$.

If the operator $T$ is an integral operator
$$
Tf(x) = \int_{\R^N} K(x,y) f(y) dm(y),
$$
then the operator $T_\mu$ is the integral operator with the same kernel 
$K$, but the integration  is performed with respect to the measure $d\mu 
= u dm$,
$$
T_\mu f(x) = \int_{\R^N} K(x,y) f(y) d\mu(y) .
$$
The adjoint of the operator $T_\mu :L^2(\mu)\to L^2(\nu)$  is the 
operator $T_\nu^*:L^2(\nu)\to L^2(\mu)$,
$$
T_\nu^* = T^* M_v.
$$
Again, in the case of integral operator, we have the following 
representation of $T^*_\nu$.
$$
T^*_\nu g(y) = \int_{\R^N} \overline{K(x,y)} g(x) d\nu(x),
$$
which explains the subscript $\nu$.

\begin{rem*}
The formula $T_\nu^* = T^* M_v$ may seem strange, it looks like the formula 
should be $M_u T^*$. However, this is a correct formula, and the 
naturally looking formula $M_u T^*$ is wrong, the main reason being 
that we are considering operators acting between different spaces.

Namely, if we represent $T_\mu$ as
$$
T_\mu = M_{v^{-1/2}} [M_{v^{1/2}} T M_{u^{1/2}} ] M_{u^{1/2}}
$$
then in brackets all operators are operators in $L^2$, so multiplication 
operators are self-adjoint ones. But outside the brackets, the operators 
$M_{v^{-1/2}}$ and $M_{u^{1/2}}$ are \emph{unitary} operators, $M_{v^{-1/2}}: 
\ran M_{v^{1/2}} \subset L^2 \to L^2(\nu)$ and $M_{u^{1/2}}: L^2(\mu) \to \ran M_{u^{1/2}} \subset 
L^2$. So their adjoint are their inverses, and so
$$
(T_\mu)^* = M_{u^{-1/2}} [ M_{u^{1/2}} T^* M_{v^{1/2}}] M_v^{1/2} = T^* 
M_v = : T_\nu^*.
$$
\end{rem*}

Let us now show that the operator $T_\mu$ is indeed a well localized 
operator in the sense of Definition \ref{df1.1}. 

Consider a decomposition of the operator $T$ with respect to the non-weighted Haar basis in $L^2(m)$, 
$$
T= \sum_{R, Q\in \cD} T\ci{R,Q}, \qquad T\ci{R, Q}: H_Q\to H_R;
$$
here recall that $H_Q$ denotes the space of all non-weighted Haar functions $h\ci Q$ on the cube $Q$  (which is a subspace of $L^2(m)$ of dimension $2^N-1$). If we chose some orthonormal bases $\{h\ci{Q,k}\}_{k=1}^{2^N-1}$ and $\{h\ci{R,j}\}_{j=1}^{2^N-1}$ in $H_Q$ and $H_R$ respectively, then we can represent $T_{R, Q}$ as 
$$
T\ci{R, Q} = \sum_{j, k=1}^{2^N-1} \La T h\ci{Q,k}, h\ci{R,j}\Ra \La\fdot, h\ci{Q,k}\Ra h\ci{R, j}.
$$
Since the rank one operator $\La\fdot, h\ci{Q,k}\Ra h\ci{R, j}$ is an integral operator with kernel $ \overline{h\ci{Q,k}(y)} h\ci{R, j}(x)$, we can conclude that $T_{R, Q}$ can be represented as an integral operator with bounded compactly supported kernel. So, the operators $T_{R, Q}$ are well defined on $L^1_\text{loc}$ functions.

Assume for a moment that $u\in L^2_\text{loc}$. Then $u\chi\ci Q\in L^2$, so 
$$
Tu\chi\ci Q = \sum_{Q',R'} T\ci{R',Q'} u\chi\ci Q
$$
where the series converges in $L^2$. Consider a cube $R$, $\ell(R)\le \ell(Q)$. If we also assume that $v\in L^2_\text{loc}$, then $\La T u\chi\ci Q , h\ci R^\nu\Ra_\nu = \La T u\chi\ci Q , v h\ci R^\nu\Ra$ is well defined and 
\begin{equation}
\label{1.3}
\La T u\chi\ci Q, h^\nu\ci R \Ra_\nu = \sum_{Q',R'\in\cD} \La T\ci{R',Q'} \chi\ci Q u, h^\nu\ci R \Ra_\nu
\end{equation}
Note, that $T\ci{R',Q'} \chi\ci Q u \ne 0$ only if $Q'\cap Q\ne \varnothing$, so in the above sum we need to consider only such $Q'$. 
Since a 
weighted Haar function $h_R^\nu$ is orthogonal in $L^2(\nu) $ to
constants, $\La T\ci{R', Q'} \chi\ci Q u, h^\nu \Ra_\nu =0$ for $R\subsetneqq R'$. Also, trivially $\La T\ci{R', Q'} \chi\ci Q u, h^\nu \Ra_\nu =0$ if $R'\cap R\ne \varnothing$. And finally, the band structure of $T$ means that $T\ci{R', Q'}=0$ if $d_\text{tree}(R',Q')>r$. So in the above sum we can consider only $R'$, $Q'$ satisfying $Q'\cap Q\ne \varnothing$, $R'\subset R$ and $d_\text{tree}(R',Q')\le r$.

Since $R'\subset R$ and $\ell(R)\le \ell(Q)$, the assumption 
$\operatorname{d}\ti{tree}(R',Q')\le r$ implies that $\ell(Q')\le 2^r 
\ell(Q)$. This together with $Q\cap Q'\ne \varnothing$ implies that 
$Q'\subset Q^{(r)}$. Indeed, $Q\cap Q'\ne \varnothing$ means that either 
$Q'\subset Q$ (and so $Q'\subset Q^{(r)}$), or $Q\subset Q'$, and in the 
latter case the estimate on the size of $Q'$ implies $Q'\subset Q^{(r)}$.

Note that if $R\not \subset Q^{(r)}$ then $R'\not \subset Q^{(r)}$ 
(because $R'\subset R$). Therefore, $\operatorname{d}\ti{tree}(R',Q')\ge 
r+2$, because we need at least one step to go above $Q^{(r)}$ and then 
at least $r+1$ steps to go to cubes of size $\ell(R)$. Therefore, $\La 
T(w^{-1}\chi\ci{Q}), h^\nu_{R}\Ra_\nu =0$ if $R\not \subset Q^{(r)}$.

A similar reasoning works for the case $\ell(R)\le 2^{-r}\ell(Q)$. 
Namely, the inclusion $R'\subset R$ and the inequality 
$\operatorname{d}\ti{tree}(R',Q')\le r$ imply that $\ell(Q')\le 
\ell(Q)$, so it follows from $Q\cap Q'\ne \varnothing$ that $Q'\subset 
Q$. If $R\not \subset Q$ then $R'\not\subset Q$ and 
$\operatorname{d}\ti{tree}(R',Q')\ge r+2$ (we need at least one step to 
go from $Q'$ above $Q$, and then at least $r+1$ steps to go to the cubes 
of size $2^{-r}\ell(R)$). Therefore, $\La T(w^{-1}\chi\ci{Q}), h_R^\nu\Ra_\nu  = 0$ 
if $\ell(R) \le 2^{-r}\ell(Q)$ and $R\not \subset Q$. 

So, the operator $T_\mu$ is \emph{lower triangularly localized}, and the same reasoning can be applied to the adjoint operator $T^*_\nu$. So we have shown, that under the assumption $u, v\in L^2_{\text{loc}}$ the operator $T_\mu$ obtained from the band operator $T$ is well localized.  \hfill\qed

To treat the general case let us note  that for $Q\cap R=\varnothing$ the above sum \eqref{1.3} has only finitely many terms. Each operator $T\ci{R',Q'}$ is an integral operator with bounded compactly supported kernel,  so $\La T\chi\ci Q u, h\ci R^\nu\Ra_\nu$ is well defined for $u, v\in L^1_\text{loc}$. 

 Moreover, if we take increasing sequences of weights $u_n, v_k\in L^2_\text{loc}$, $u_n\nearrow u$, $v_k\nearrow v$, and define $d\nu_k =v_k dm$, then 
 $$
 \lim_{n\to\infty} \lim_{k\to \infty} \La T \chi\ci Q u_n , h\ci R^{\nu_k}\Ra_{\nu_k} = 
 \La T\chi\ci Q u, h\ci R^\nu\Ra_\nu = \La M_{v^{1/2}} T M_{u^{1/2}} \chi\ci Q u^{1/2}, h^\nu\ci R v^{1/2} \Ra,  
$$
where the function $M_{v^{1/2}} T M_{u^{1/2}} \chi\ci Q u^{1/2}$ is interpreted exactly as in the Remark \ref{r0.2}
 
\subsection{Sawyer type results for well localized operators}
The following two theorems can be also considered  the main result of the paper.

\begin{thm}
\label{t1.1}
Let $T=T_\mu$ be a well localized operator acting (formally) from 
$L^2(\mu) $ to $L^2(\nu)$. Then $T_\mu$ is a bounded operator from 
$L^2(\mu) $ to $L^2(\nu)$ if and only if $T$ and its formal adjoint 
$T^*_\nu$ are uniformly bounded on characteristic functions of cubes, 
i.e.~iff
\begin{align*}
\|T_\mu\chi\ci Q\|\ci{L^2(\nu)}^2  & \le C \|\chi\ci Q\|\ci{L^2(\mu)}^2 
= C\mu(Q),\qquad \forall Q\in \cD,  \\
\|T_\nu^*\chi\ci Q\|\ci{L^2(\mu)}^2 &\le C \|\chi\ci Q\|\ci{L^2(\nu)}^2 
= C\nu(Q), \qquad \forall Q\in \cD.
\end{align*}
Moreover, the norm of $T$ can be estimated by a constant depending  only on the dimension $N$, the above constants $C$ and $r$ from the definition of well localized operator. 
\end{thm}

Theorem \ref{t0.1} is an immediate corollary of this theorem.

The assumptions that $T_\mu$ and $T^*_\nu$  are uniformly bounded on the 
characteristic functions of cubes can be relaxed a little: one does not 
have to integrate $T_\mu\chi\ci Q$ over the whole space.
Namely, Theorem \ref{t1.1} can be restated as follows

\begin{thm}
\label{t1.2}
Let $T=T_\mu$ be a well localized operator acting (formally) from 
$L^2(\mu) $ to $L^2(\nu)$. Then $T_\mu$ is a bounded operator from 
$L^2(\mu) $ to $L^2(\nu)$ if and only if the following two conditions
hold:
\begin{enumerate}
	\item $\left|\La T_\mu \chi\ci Q, \chi\ci R \Ra_\nu \right| \le 
C\mu(Q)^{1/2}\nu(R)^{1/2}$ for all cubes $Q$, $R$ of comparable size, 
$2^{-r}\ell(Q)\le \ell(R)\le 2^r\ell(Q)$; here $r$ is the number from 
the definition of well localized operator.
	\item  For all cubes $Q$
	$$
	\int_Q\left| T_\mu \chi\ci Q \right|^2 \,d\nu \le C\mu(Q), \qquad
	\int_Q\left| T_\nu^* \chi\ci Q \right|^2 \,d\nu \le C\nu(Q).
	$$
\end{enumerate}
Moreover, the norm of $T$ can be estimated by a constant depending  only on the dimension $N$, the above constants $C$ and $r$ from the definition of well localized operator. 
\end{thm}
Theorem \ref{t0.3} is an immediate corollary of the above Theorem \ref{t1.2}.

Note that in the condition (2) of the theorem one can replace $T_\mu 
\chi\ci Q$ by its orthogonal (in $L^2(\nu)$) projection onto the 
subspace of functions $F$ with zero average over $Q$, $\int_Q 
f\,d\nu=0$, and similarly for $T_\nu^*\chi\ci Q$.  The condition (1)  implies that $\left|\La T_\mu\chi\ci Q, 
\chi\ci Q\Ra_\nu \right|\le C(\mu(Q))^{1/2}(\nu(Q))^{1/2}$, i.e. that the projection onto the 
orthogonal complement of this subspace is bounded.

While in our main example the measures $\mu$ and $\nu$ are absolutely continuous, and the operator $T$ came from an operator in the non-weighted $L^2$, in Theorems \ref{t1.1} and \ref{t1.2} the measures are arbitrary regular Borel measures and $T_\mu$ is an arbitrary operator $L^2(\mu)\to L^2(\nu)$. We only need to know its bilinear form $\La T_\mu \chi\ci Q, \chi\ci R\Ra_\nu$ to apply the theorems.

\begin{rem}
We do not have to assume that the operator $T$ is  bounded in the non-weighted $L^2$ to be able to apply Theorem \ref{t1.1} and \ref{t1.2} to get the weighted norm inequalities  for  $T$. We only need to be able to define $\La T_\mu \chi\ci Q, \chi\ci R \Ra_\nu $. In Section \ref{s1.1} above we had shown how one can treat this expression if the operator $T$ is bounded in the non-weighted $L^2$. The same reasoning will work if $T$ is bounded in any non-weighted $L^p$.

In general situation (for example, if we have unbounded Haar multipliers) some care is needed to interpret $\La T_\mu \chi\ci Q , \chi\ci R\Ra_\nu$, but after that one can apply Theorem \ref{t1.1} and \ref{t1.2}. 
\end{rem}
\section{Proof of Theorems \ref{t1.1} and \ref{t1.2}.}
\label{s2}

\subsection{Martingale difference decomposition}

Denote by $E_{k}=E_k^\mu$ the averaging operator in $L^2(\mu)$ over 
dyadic cubes of  size
(length of the side) $2^{k}$, namely $E_{k}^\mu f(x) = \mu(Q)^{-1}
\int_{Q} f d\mu$, where $Q$ is a dyadic cube of size $2^{k}$
containing
$x$ (for the sake of definiteness, we consider cubes of  the form
$x_0+[a,b)^N$). If $Q$ is a cube  of size $2^{k}$, we  denote by
$E\ci Q^\mu
f$ the restriction of $E_{k}^\mu f$ to $Q$: $E\ci Q^\mu  f := (\mu(Q)^{-1}
\int_{Q} f d\mu) \chi\ci Q = \chi\ci Q E_{k}^\mu f $.

Let $\Delta_{k}=\Delta_k^\mu := E_{k-1}^\mu-E_{k}^\mu$.
Again for a  dyadic cube
$Q$ of size $2^{k}$,
denote by $\Delta_{Q}^\mu f$ the restriction of $\Delta_{k}^\mu f $ to
$Q$. Clearly, for any $f\in L^{2}(\mu)$, the functions $ \Delta_{Q}^\mu
f$, $Q\in\mathcal{D}$, are orthogonal to each other, and that for any
fixed
$n$ we have the orthogonal decomposition
\begin{align}
\label{3.1}
f & = \sum_{Q\in\cD, \ell(Q)\le 2^{n}} \Delta_{Q}^\mu f +
\sum_{Q\in\cD, \ell(Q)= 2^{n}}  E_{Q}^\mu f ,
\\
\notag
\|f\|^{2}_{L^{2}(\mu)} &= \sum_{Q\in\cD, \ell(Q)\le 2^{n}} \|
\Delta_{Q}^\mu f \|^{2} + \sum_{Q\in\cD, \ell(Q)= 2^{n}} \|
E_{Q}^\mu f \|^{2}.
\end{align}

\subsection{Paraproducts}

Define the \emph{paraproduct} $\Pi^\mu=\Pi^\mu_T$, acting (formally) 
from  $L^2(\mu)$ to  $L^2(\nu)$ by
$$
\Pi^\mu f := \sum_{Q\in \cD} E_Q^\mu f  \sum\begin{Sb}R\in \cD,\ 
R\subset Q, \\ \ell(R)=2^{-r}\ell(Q)\end{Sb} \Delta_R^\nu T_\mu \chi\ci Q.
$$
The paraproduct $\Pi^\nu=\Pi^\nu_{T^*}$ is defined similarly,
$$
\Pi^\nu f := \sum_{Q\in \cD} E_Q^\nu f  \sum\begin{Sb}R\in \cD,\ 
R\subset Q, \\ \ell(R)=2^{-r}\ell(Q)\end{Sb} \Delta_R^\mu T_\nu^* \chi\ci
Q.
$$

\begin{rem}
\label{r2.1}
Note, that it follows from the definition of well localized operator 
that if $R\subset Q$, and $\ell(R) \le 2^{-r} \ell(Q)$, then for any 
$Q'\supset Q$
$$
\Delta_R^\nu T_\mu \chi\ci{Q'} = \Delta_R^\nu T_\mu \chi\ci{Q}.
$$
In other words, in the definition of $\Pi$ we can always replace 
$\chi\ci Q$ by $\chi\ci{Q'}$ where $Q'$ is a bigger cube.

That essentially means that formally we can write $T1$ instead of 
$T\chi\ci  Q$, so the definition is more in line with the standard 
definition of a paraproduct.
\end{rem}

The following lemma describes  the matrix of $\Pi$ with respect to the 
weight\-ed Haar systems in $L^2(\mu)$ and $L^2(\nu)$.

\begin{lm}
\label{l2.2}
Let $Q$, $R$ be dyadic cubes. Then for the paraproduct $\Pi=\Pi^\mu$ 
defined above
\begin{enumerate}
	\item If $\ell(R) \ge  2^{-r}\ell(Q)$ then $\La \Pi^\mu h^\mu_Q, 
h^\nu_R\Ra_\nu = 0$ for all weighted Haar functions $h^\mu_Q$ and 
$h^\nu_R$.

\item If $R\not\subset Q$, then 
$\La \Pi^\mu h^\mu_Q, h^\nu_R\Ra_\nu = 0$ for all weighted Haar 
functions $h^\mu_Q$ and $h^\nu_R$.

\item If $\ell(R) <2^{-r} \ell(Q)$, then for all weighted Haar functions 
$h^\mu_Q$ and $h^\nu_R$
$$
\La \Pi^\mu h_Q^\mu, h_R^\nu\Ra_\nu = \La T h_Q^\mu, h_R^\nu\Ra_\nu;
$$
in particular, if $R\not\subset Q$, then both sides of the equality are 
$0$.

\end{enumerate}
\end{lm}

\begin{proof}
Let us use $Q'$ and $R'$ for the summation indices in the paraproduct, 
i.e.~let us write
$$
\Pi^\mu h_Q^\mu := \sum_{Q'\in \cD} E_{Q'}^\mu h_Q^\mu 
\sum\begin{Sb}R'\in \cD,\ R'\subset Q', \\ 
\ell(R')=2^{-r}\ell(Q')\end{Sb} \Delta_{R'}^\nu T_\mu \chi\ci{Q'}.
$$
Since $h^\nu_R$ is orthogonal to ranges of all projections 
$\Delta_{R'}^\nu$ except $\Delta_R^\nu$ we can write
\begin{equation}
\label{2.1}
\La \Pi h^\mu_Q, h_R^\nu\Ra_\nu = \La (E_{Q'}^\mu h^\mu_Q) \Delta_R^\nu 
T\chi\ci{Q'}, h_R^\nu \Ra_\nu = a\La T\chi\ci{Q'}, h_R^\nu \Ra_\nu 
\end{equation}
where $Q'$ is the grandparent of $R$ of order $r$ (i.e.~the cube $Q'$, 
$Q'\supset R$ and such that $\ell(Q')= 2^r \ell(R)$) and $a$ is the value of $E_{Q'} h^\mu_{Q}$ on $Q'$, $E_{Q'} h^\mu_{Q} = a \chi\ci{Q'}$.

It is easy to see that $E^\mu_{Q'} h^\mu_Q \not\equiv 0$ (equivalently $a\ne 0$) only if $Q'\subsetneqq 
Q$.  
Therefore, see \eqref{2.1},  $\La \Pi h^\mu_Q, h_R^\nu\Ra_\nu\ne )$ only if $Q'\subsetneqq Q$
 and statements 1 and 2 of the lemma follow immediately. 

Indeed, if $\ell(R)\ge 2^{-r}\ell(Q)$ and $\ell(Q') =2^r \ell(R)$, the inclusion $Q'\subsetneqq 
Q$ is impossible, so  $ \La \Pi h^\mu_Q, h_R^\nu\Ra_\nu = 0$, and the statement (1) is proved. 

If  $R\not\subset Q$, then the inclusion $Q'\subsetneqq 
Q$ implies that $R\not\subset Q'$. But $\ell(R) = 2^{-r}\ell(Q')$, so by the definition of a well localized operator $ \La \Pi h^\mu_Q, h_R^\nu\Ra_\nu = 0$ and the statement (2) is proved.

Let us prove statement 3. Let $\ell(R)<2^{-r}\ell(Q)$. If $R\not\subset 
Q$ then by the statement 2 of the lemma $\La \Pi h^\mu_Q, h^\nu_R\Ra_\nu 
=0$, and $\La T_\mu h^\mu_Q, h^\nu_R\Ra_\nu =0$ by the definition of 
well localized operators (``children'' $Q_k$ of $Q$ are cubes of size at 
least $2^r\ell(R)$, and it follows from the definition of well localized 
operator that $\La T_\mu \chi\ci{Q_k}, h^\nu_R\Ra_\nu =0$). So, we only 
need to consider the case $R\subset Q$.

Let $Q_1$ be the ``child'' of $Q$ containing $R$ (i.e.~$R\subset 
Q_1\subset Q$, $\ell(Q_1) =\ell(Q)/2$), and let $a$ be the value of 
$h^\mu_Q$ on $Q_1$. Then, since for all other children $Q_k$ of $Q$ the 
definition of well localized operator implies $\La T\chi\ci{Q_k} , 
h^\nu_R\Ra_\nu =0$, we can conclude that
$$
\La T h_{Q}^\mu , h^\nu_R\Ra_\nu = a \La T\chi\ci{Q_1} , h^\nu_R\Ra_\nu
$$
On the other hand we have shown before, see \eqref{2.1} that
$$
\La \Pi h^\mu_Q, h_R^\nu\Ra_\nu = \La (E_{Q'}^\mu h^\mu_Q) \Delta_R^\nu 
T\chi\ci{Q'}, h_R^\nu \Ra_\nu
$$
where $R\subset Q' \subsetneqq Q$, $\ell(Q')=2^{r}\ell(R)$. Therefore 
$Q'\subset Q_1$ and so $E_{Q'}^\mu h_Q^\mu =a\chi\ci{Q'}$. We also know, see Remark 
\ref{r2.1}, that because $Q'\subset Q_1$ we have equality $\Delta_R^\nu 
T\chi\ci{Q'} = \Delta_R^\nu T\chi\ci{Q_1}$. Thus we can continue:
$$
\La \Pi h^\mu_Q, h_R^\nu\Ra_\nu  =
a\La  \Delta_R^\nu T\chi\ci{Q'}, h_R^\nu \Ra_\nu = a\La  \Delta_R^\nu 
T\chi\ci{Q_1}, h_R^\nu \Ra_\nu= a\La  
T\chi\ci{Q_1}, h_R^\nu \Ra_\nu.
$$
Therefore $\La \Pi h^\mu_Q, h_R^\nu\Ra_\nu = \La T h_{Q}^\mu , 
h^\nu_R\Ra_\nu$, and the lemma is proved.
\end{proof}  

It follows from Lemma \ref{l2.2} and the definition of a well localized 
operator, that in the weighted Haar bases the matrix of the difference $T_\mu - \Pi^\mu_T - (\Pi^\nu_{T^*})^*$ 
has finitely many diagonals, which are bounded by the assumption (1) of 
Theorem \ref{t1.2}. Therefore, the main part in the proof Theorem \ref{t1.2} (and thus 
Theorem \ref{t1.1}) is  to prove that the paraproducts 
are bounded. Of course, one also need to take care of the terms like $E\ci Q^\mu f$ in the decomposition \eqref{3.1}, but this is ratther simple, we present the details in Section \ref{s3.4} below

\subsection{Boundedness of the paraproduct}

We will need the following well known theorem.

Let $f\ci{\!R}:= \frac1{\mu(R)} \int_R f\,d\mu$ be the average of the 
function $f$ with respect to the measure $\mu$.

\begin{thm}[Dyadic Carleson Embedding Theorem]
\label{t2.3}
If the numbers $a\ci Q\ge 0$, $Q\in \cD$ satisfy the following Carleson 
measure condition
\begin{equation}
\label{2.2}
\sum_{Q\subset R} a\ci Q \le \mu(R),
\end{equation}
then for any $f\in L^2(\mu)$
$$
\sum_{R\in \cD} a\ci R |f\ci{\!R}|^2 \le C \|f\|_{L^2(\mu)}^2
$$
where $C$ is an absolute constant.  
\end{thm}
This theorem is very well known, cf \cite{Duren1970}. For an alternative prove see also \cite{NTV}, 
\cite{NT}, \cite{NTV-name}, where it was proved with the constant $C=4$ 
using Bellman function method. It was also proved in \cite{NTV-name} 
that the constant $C=4$ is optimal. We should mention, that in 
\cite{NT}, \cite{NTV-name} this theorem was proved for $\R^1$, but the 
same proof works for general martingale setup. A proof for $\R^2$ was 
presented in \cite{NTV}, and the same proof works for $\R^N$.

Let us now show that the paraproduct $\Pi=\Pi^\mu_T$ is bounded. Ranges 
of the projections $\delta^\nu_R$ are mutually orthogonal, so to prove 
the boundedness of the paraproduct $\Pi^\mu_T$ it is sufficient to show 
that the numbers
$$
a\ci Q := \sum\begin{Sb}R\in \cD, R\subset Q\\ \ell(R) = 
2^{-r}\ell(Q)\end{Sb} \|\Delta^\nu_R T_\mu \chi\ci R \|^2_{L^2(\nu)} \le
$$
satisfy the Carleson Measure Condition  \eqref{2.2} from Theorem 
\ref{t2.3}. Let us prove this.

Consider a cube $\wt Q$. We want to show that
$$
\sum_{Q\subset \wt Q} \sum\begin{Sb}R\in \cD, R\subset Q\\ \ell(R) = 
2^{-r}\ell(Q)\end{Sb} \|\Delta^\nu_R T_\mu\chi\ci Q \|_{L^2(\nu)}^2 \le 
C \mu(Q).
$$
By Remark \ref{r2.1} we can replace $\chi\ci Q$ by $\chi\ci{\wt Q}$, so 
the desired estimates becomes
$$
\sum\begin{Sb}R\in \cD, R\subset \wt Q\\ \ell(R) \le 2^{-r}\ell(\wt 
Q)\end{Sb} \|\Delta^\nu_R T_\mu\chi\ci{\wt Q} \|_{L^2(\nu)}^2 \le 
\sum_{R\subset \wt Q} \|\Delta^\nu_R T_\mu\chi\ci{\wt Q} \|_{L^2(\nu)}^2
\le \|\chi\ci{\wt Q} T_\mu\chi\ci{\wt Q} \|_{L^2(\nu)}^2 .
$$
By assumption (2) of Theorem \ref{t1.2}
$$
\|\chi\ci{\wt Q} T_\mu\chi\ci{\wt Q} \|_{L^2(\nu)}^2 :=\int_{\wt Q} | 
T_\mu\chi\ci{\wt Q} |^2d\nu \le C\mu(\wt Q)
$$
and so the sequence $a\ci Q$, $Q\in \cD$ satisfies the condition 
\eqref{2.2} \hfill\qed

\subsection{Why $T$ is bounded.}
\label{s3.4}
Let $f\in L^2(\mu)$, $g\in L^2(\nu)$.

We want to estimate
$$|\La T_\mu f, g\Ra_\nu| \le C\|f\|_{L^2(\mu)}\|g\|_{L^2(\nu)}.
$$
It is sufficient to prove the estimate on a dense set of compactly 
supported functions. Each compact subset of $\R^N$ is contained in at 
most $2^N$ cubes of the same size as the size of this compact subset, so let $Q_k$, $k=1, 2, \ldots, 2^N$ 
be cubes of size $2^d$ containing supports of $f$ and $g$.

Let us decompose $f$ and $g$ using the martingale difference
decomposition,
$$
f= \sum_{k=1}^{2^N} E_{Q_k}^\mu f + \sum_{Q:\ell(Q)\le 2^d}\Delta^\mu_Q 
f, \qquad
g= \sum_{k=1}^{2^N} E_{Q_k}^\nu g + \sum_{Q:\ell(Q)\le 2^d}\Delta^\nu_Q g
.
$$
Note, that $\Delta_Q^\mu f$ and $\Delta^\nu_Q g$ are $\mu$- and 
$\nu$-Haar functions on the cube $Q$.

By Lemma \ref{l2.2}
\begin{multline*}
\Bigl\La T_\mu \sum_{Q} \Delta^\mu_Q f, \sum_R \Delta_R^\nu g 
\Bigr\Ra_\nu =
\La \Pi_T^\mu f, g \Ra_\nu + \La f, \Pi_{T^*}^\nu g \Ra_\nu \\ + 
\sum_{Q:\ell(Q)\le 2^d} \left( \sum_{R: 2^{-r}\ell(Q) \le \ell(R) \le 
2^r \ell(Q)}  \La T_\mu \Delta_Q^\mu f, \Delta_R^\nu g\Ra_\nu \right)
\end{multline*}
We know that the paraproducts $\Pi_T^\mu$ and $\Pi_{T^*}^\nu$ are 
bounded, so it remains to estimate the last sum. It follows from the 
assumption (1) of Theorem \ref{t1.2} that $| \La T_\mu \Delta_Q^\mu f, 
\Delta_R^\nu g\Ra_\nu | \le C \|\Delta_Q^\mu f\|_{L^2(\mu)} \,\| \Delta_R^\nu 
g\|_{L^2(\nu)}$ if $2^{-r}\ell(Q) \le \ell(R) \le 2^r \ell(Q)$.

On the other hand, it follows from the definition of a well localized 
operator, that given a cube $Q$ at most $M$ terms $ \La T_\mu 
\Delta_Q^\mu f, \Delta_R^\nu g\Ra_\nu$ are non-zero, where $M$ does not 
depend on the choice of $Q$. Therefore we can split the sum into at most 
$M$ sums of form (operators with ``one diagonal'')
$$
\sum_Q \La T_\mu \Delta_Q^\mu f, \Delta_{R(Q)}^\nu g\Ra_\nu,
$$
which can be estimated by Cauchy--Schwarz inequality.

It remains to estimate the terms involving $E_Q$. By assumption (1) of 
Theorem \ref{t1.2} $|\La T_\mu E_{Q_k}^\mu f , E_{Q_j}^\nu g \Ra | \le C 
\|\chi\ci{Q_k}\|_{L^2(\mu)} \|\chi\ci{Q_j}\|_{L^2(\nu)} \le C\|f\|_{L^2(\mu)} \|g\|_{L^2(\nu)}$, and 
since there are at most $2^N$ cubes $Q_k$ we can estimate
$$
\sum_{k,j} |\La  T_\mu E_{Q_k}^\mu  f , E_{Q_j}^\nu g \Ra_\nu  |\le 
2^{2N} C \|f\|_{L^2(\mu)} \|g\|_{L^2(\nu)} .
$$

Let us now estimate
$$
\sum_{k} \Bigl|\sum_{R} \La T_\mu E_{Q_k}^\mu f, \Delta_R^\nu g \Ra_\nu \Bigr|
$$
For a fixed $Q_k$ we have by the assumption of Theorem \ref{t1.1}
$$
\Bigl|\sum_{R: \ell(R) < \ell(Q_k)}  \La T_\mu E_{Q_k}^\mu f, 
\Delta_R^\nu g \Ra_\nu \Bigr| \le C \|f\|_{L^2(\mu)} \|g\|_{L^2(\nu)}.
$$
But we have at most $2^N$ cubes 
$Q_k$, so the whole sum is bounded by $2^N C\|f\|_{L^2(\mu)} \|g\|_{L^2(\nu)}$. 

The sum
$$
\sum_{k} \Bigl|\sum_{R} \La T_\mu \Delta_{R}^\mu f, E_{Q_k}^\nu g \Ra_\nu \Bigr|
$$
can be estimated similarly. 

\hfill\qed

\end{document}